\documentclass[12pt]{amsart}
\usepackage{amsfonts}
\usepackage{amssymb}
\theoremstyle{definition}

\theoremstyle{remark}

\newcommand{\X}{\mathfrak X}

\begin{document}

\centerline{\large\bf ON THE SECTIONAL CURVATURE}
\centerline{\large\bf OF K\"AHLER MANIFOLDS OF INDEFINITE METRICS
\footnote{\it PLISKA Studia mathematica bulgarica. Vol. 9, 1987, p. 66-70}}

\vspace{0.2in}
\centerline{\large OGNIAN  KASSABOV, ADRIJAN BORISOV}

\vspace{0.3in}
{\sl  In this paper we study some curvature properties of K\"ahler manifolds of
indefinite metrics.}

\vspace{0.2in}
{\bf 1. Introduction.} K u l k a r n i [5] has proved, that if a connected manifold of 
indefinite metric has bounded from above or from below sectional curvature, then it is
of constant sectional curvature. Similar results for indefinite K\"ahlerian metrics
are found in [1].

H a r r i s [4] and, independently, D a j c z e r and N o m i z u [3] generalise [5] by
using restrictions only on the sectional cuvatures for the timelike planes (or for the 
spacelike planes). The purpose of this paper is to prove analogous results for K\"ahler
manifolds of infedinite metrics.

Let $M$ be a K\"ahler manifold of indefinite metric tensor $g$ and complex structure $J$.
Let $\nabla$ denote the covariant differentiation with respect to the Levi-Civita
connection. Then $g(JX,JY)=g(X,Y)$ for $X,Y \in \X(M)$, $\nabla J=0$ and the curvature tensor 
$R$ \,has the properties:

1) $R(X,Y,Z,U)=-R(Y,X,Z,U)$;

2) $R(X,Y,Z,U)+R(Y,Z,X,U)+R(Z,X,Y,U)=0$;

3) $R(X,Y,Z,U)=-R(X,Y,U,Z)$;

4) $R(X,Y,Z,U)=R(X,Y,JZ,JU)$.

A pair $\{ x,y \}$ of tangent vectors at a point $p\in M$ is said to be orthonormal of
signature $(+,-)$ if $g(x,y)=0$, $g(x,x)=1$, $g(y,y)=-1$. Analogously are defined pairs of
signature (+,+) or $(-,-)$. Similarly, one speaks of the signature of a 2-plane $\sigma$,
depending on the signature of the restriction of the metric on $\sigma$ [2].

We recall that the curvature of a nondegenerate 2-plane $\sigma$ with a basis $\{ x,y \}$
is defined by 
$$
	K(\sigma) = \frac{R(x,y,y,x)}{\pi_1(x,y,y,x)} \ ,
$$
where
$$
	\pi_1(x,y,z,u)=g(x,u)g(y,z)-g(x,z)g(y,u) \ .
$$

A plane $\sigma$ is said to be holomorphic (resp. antiholomorphic) if $\sigma=J\sigma$
(resp. $\sigma \perp J\sigma$). Then a 2-plane $\sigma$  is holomorphic (resp. antiholomorphic)
if and only if it has a basis $\{ x,Jx\}$ (resp. $\{ x,y\}$ with $g(x,y)=g(x,Jy)=0$). We call the 
pair $\{ x,y \}$ antihilomorphic of signature $(+,-)$, if $g(x,x)>0$, $g(y,y)<0$ and
$g(x,y)=g(x,Jy)=0$. Similarly we can define antiholomorphic pairs of signature (+,+) 
or $(-,-)$.  

If $\sigma$ is a nondegenerate 2-plane with an orthonormal basis  $\{ x,y \}$  we denote
$$
	K(x,y)=K(\sigma)
$$
and in particular, if $x$ is a unit vector, i.e. $g(x,x)=1$ or $g(x,x)=-1$ we denote
$$
	H(x)=K(x,Jx) \ .
$$

The manifold $M$ is said to be of constant holomorphic sectional curvature if for
any point $p \in M$ the curvature of an arbitrary holomorphic 2-plane $\sigma$ in
$T_p(M)$ doesn't depend on $\sigma$. As in the case of definite metric this occurs 
when and only when the curvature tensor $R$ has the form
$$
	R=\frac{\mu}{4}(\pi_1+\pi_2) \ ,
$$
where
$$
	\pi_2(x,y,z,u)=g(x,Ju)g(y,Jz)-g(x,Jz)g(y,Ju)-2g(x,Jy)g(z,Ju) \ .
$$
In this case $\mu$ is a global constant if $M$ is connected.

\vspace{0.2in}
{\bf 2. Preliminary considerations.} All the manifolds under consideration in this 
and the following section will be assumed to be connected.

\vspace{0.1in}
P\,r\,o\,p\,o\,s\,i\,t\,i\,o\,n 1. {\it Let $M$ be a $2m$-dimensional K\"ahler 
manifold of indefinite metric, $m\ge 2$. If
$$
	R(x,Jx,Jx,y)=0 \ ,   \leqno (1)
$$
whenever the pair $\{ x,y \}$ is antiholomorphic of signature $(+,-)$, then $M$
is of constant holomorphic sectional curvature. }

\vspace{0.1in}
P\,r\,o\,o\,f. Let $p\in M$ and $x,y$ be unit vectors in $T_p(M)$, such that 
$\{ x,y \}$ is an antiholomorphic pair of signature $(+,-)$. According to (1),
for any $\alpha$, $|\alpha|<1$
$$
	R(x+\alpha y, Jx+\alpha Jy,Jx+\alpha Jy, \alpha x+y)=0
$$
holds good and hence using (1), we obtain
$$
	\alpha\{ H(x)-K(x,y)-3K(x,Jy) \}
$$
$$
	+\alpha^3\{ H(y)-K(x,y)-3K(x,Jy)\}+(3+\alpha^2)\alpha^2R(y,Jy,Jy,x)=0 \ .
$$
Hence we derive
$$
	R(y,Jy,Jy,x)=0 \ , \leqno (2)
$$
$$
	 H(x)-K(x,y)-3K(x,Jy)+\alpha^2\{  H(y)-K(x,y)-3K(x,Jy) \} =0  \leqno (3)
$$
for $|\alpha|<1$, $\alpha \ne 0$ and hence (3) holds for any $\alpha$, $|\alpha|<1$. 
Now let $\alpha =0$. Then
$$
  H(x)=K(x,y)+3K(x,Jy) \ .
$$
Hence we find
$$
	K(x,y)=K(x,Jy) \ ,  \leqno (4)
$$
$$
	H(x)=4K(x,y) \ .    \leqno (5)
$$
From (3), (4), (5) it follows
$$
	H(x)=H(y) \ .   \leqno (6)
$$

If $m=2$ we put $\mu = H(x)$. Then from (1), (2), (4), (5), (6) and the properties of
the curvature tensor we obtain $R=(\mu/4)(\pi_1+\pi_2)$, which proves the assertion.

Let $m>2$. We choose a unit vector $z \in T_p(M)$, such that span$\{x,y,z\}$ is 
antiholo\-morphic. Then (6) implies $H(x)=H(y)=H(z)$, i.e. (6) holds also if the unit vectors 
$x,y$ form an antiholomorphic pair of signature (+,+) or $(-,-)$. Now let $\{u,v\}$ be
arbitrary unit vectors in $T_p(M)$. We choose a unit vector $x$ in $T_p(M)$ such that
span$\{x,u\}$ and span$\{x,v\}$ are antiholomorphic. Then (6) implies
$$
	H(u)=H(x)=H(v) \ ,
$$ 
thus proving Proposition 1.

\vspace{0.1in}
R\,e\,m\,a\,r\,k 1. Obviously, the signature of $\{x,y\}$ in Proposition 1 can be replaced by
$(-,+)$.

\vspace{0.1in}
P\,r\,o\,p\,o\,s\,i\,t\,i\,o\,n 2. {\it Let $M$ be a $2m$-dimensional K\"ahler manifold of
indefinite metric, $m>2$. If the signature of $M$ is $(-,-,...,+,+,+,+,...)$ and (1) holds
whenever the pair $\{x,y\}$ is antiholomorphic of signature $(+,+)$, then $M$ is of
constant holomorphic sectional curvature.}

\vspace{0.1in}
P\,r\,o\,o\,f. In a point $p \in M$ let $\{x,y\}$ be an arbitrary antiholomorphic pair of
signature $(+,-)$. We choose a vector $z$ in $T_p(M)$, such that $\{x,z\}$ is an
antiholomorphic pair of signature $(+,+)$. Then for sufficiently large $\alpha$ the pair
$\{x,y+\alpha z\}$ is antiholomorphic of signature $(+,+)$. Now from
$$
	R(x,Jx,Jx,y+\alpha z)=0 \ ,
$$
we obtain the conditions of Proposition 1. Hence Proposition 2 follows.

\vspace{0.1in}
R\,e\,m\,a\,r\,k 2. If the signature of $M$ is $(-,-,-,-,...,+,+,...)$, the signature of
$\{x,y\}$ in Proposition 2 can be replaced by $(-,-)$.

\vspace{0.1in}
P\,r\,o\,p\,o\,s\,i\,t\,i\,o\,n 3. {\it Let $M$ be a $2m$-dimensional K\"ahler manifold of
indefinite metric, $m>2$. If
$$
	R(x,y,y,z)=0
$$
whenever \, {\rm span}$\{ x,y,z \}$ is antiholomorphic and the pair $\{x,y\}$ is of signature $(+,-)$, \
then $M$ is of constant holomorphic sectional curvature. }

\vspace{0.1in}
P\,r\,o\,o\,f. Let e.g. the signature of $M$ be $(-,-,-,-,...,+,+,...)$ and let $\{x,y\}$
be an arbitrary antiholomorphic pair of signature $(+,-)$ in a point $p\in M$. We choose 
$z$ in $T_p(M)$ such that span$\{x,y,z\}$ is antiholomorphic and $g(z,z)<0$. Then we 
can choose $\alpha \ne 0$ such that the antiholomorphic pair $\{x+\alpha z,\alpha Jx+Jz\}$
is of signature $(+,-)$. Then 
$$
	R(x+\alpha z,\alpha Jx+Jz,\alpha Jx+Jz,y)=0   \leqno (7)
$$
and since (7) holds also if we change $\alpha$ by $-\alpha$, we obtain
$$
	R(x,Jx,Jx,y)=0 \ .
$$
According to Proposition 1, $M$ is of constant holomorphic sectional curvature.

\vspace{0.1in}
R\,e\,m\,a\,r\,k 3. The signature of the pair $\{x,y\}$ in Proposition 3 can be replaced by 
$(-,+)$ or by $(+,+)$ or $(-,-)$ in the case of an appropriate signature of $M$.

\vspace{0.1in}
R\,e\,m\,a\,r\,k 4. We can define in an obvious manner triples of signature $(+,+,+)$,
$(+,+,-)$, etc. Then in the case of an appropriate signature of $M$ Proposition 3 can
be formulated for antiholomorphic triples $\{x,y,z\}$ of signature $(+,+,+)$,
$(+,+,-)$, etc.

\vspace{0.1in}
{\bf 3. Main results.}

\vspace{0.1in}
T\,h\,e\,o\,r\,e\,m 1. {\it Let $M$ be a K\"ahler manifold of indefinite metric, $m\ge 2$.
If for each point $p\in M$ there exists a constant $c(p)$, such that for any vector
$x\in T_p(M)$ with $g(x,x)=1$ the holomorphic sectional curvature $H(x)$ satisfies
$|H(x)| \le c(p)$, then $M$ is of constant holomorphic sectional curvature.}

\vspace{0.1in}
P\,r\,o\,o\,f. Let $p\in M$ and $\{x,y\}$ be arbitrary orthonormal pair of signature
$(+,-)$ in $p$, such that $g(x,y)=g(x,Jy)=0$. Then from
$$
	| H( \frac{x+\alpha y}{\sqrt{1-\alpha^2}} ) | \le c(p)
$$
for $|\alpha|<1$ we obtain 
$$
	\begin{array}{c}
		|H(x)+4\alpha R(x,Jx,Jx,y)-2\alpha^2\{ K(x,y)+3K(x,Jy)\}  \\
		+4\alpha^3R(x,Jy,Jy,y)+\alpha^4H(y)| \le c(p)(1-\alpha^2)^2 \ .
	\end{array} \leqno (8)
$$
Hence we obtain easily
$$
	|H(x)-2\alpha^2\{ K(x,y)+3K(x,Jy)\}+\alpha^4H(y)| \le c(p)(1-\alpha^2)^2  \ .   \leqno (9)
$$
On the other hand from (8) it follows by continuity
$$
	\begin{array}{c}
		H(x)+4\alpha R(x,Jx,Jx,y)-2\alpha^2\{ K(x,y)+3K(x,Jy)\}  \\
		+4\alpha^3R(x,Jy,Jy,y)+\alpha^4H(y) =0 
	\end{array} 
$$
for $\alpha =\pm 1$, which implies
$$
	R(x,Jx,Jx,y)+R(x,Jy,Jy,y)=0 \ .  \leqno (10)
$$
From (8) and (10) we find
$$
	\begin{array}{c}
		|H(x)-2\alpha^2\{ K(x,y)+3K(x,Jy)\} +4\alpha(1-\alpha^2)R(x,Jx,Jx,y) \\
		+\alpha^4H(y)| \le c(p)(1-\alpha^2)^2 
	\end{array} 
$$
and hence and (9)
$$
	|\alpha R(x,Jx,Jx,y)| \le \frac 12 c(p)(1-\alpha^2) \ .
$$
Now by continuity we obtain \ $R(x,Jx,Jx,y)=0$ \ and the Theorem follows from Propo\-sition 1.

\vspace{0.1in}
R\,e\,m\,a\,r\,k 5. In Theorem 1 the requirement $g(x,x)=1$ can be changed by $g(x,x)=-1$.

\vspace{0.1in}
T\,h\,e\,o\,r\,e\,m 2. {\it Let $M$ be a K\"ahler manifold of indefinite metric, $m> 2$.
If for each point $p\in M$ there exists a constant $c(p)$, such that for any antiholomorphic
2-plane $\sigma$ in $ T_p(M)$ \ $K(\sigma)\le c(p)$ holds,
then $M$ is of constant holomorphic sectional curvature.}

\vspace{0.1in}
P\,r\,o\,o\,f. Let $p \in M$ and $\{x,y\}$ be an arbitrary orthonormal pair of signature 
$(+,-)$ and $g(x,Jy)=0$. We choose a unit vector $z \in T_p(M)$ such that span$\{x,y,z\}$
is anti\-holomorphic. Then
$$
	K({\rm span}\{x+\alpha y,z\}) \le c(p)
$$
for $\alpha \ne \pm 1$ implies
$$
	K(x,z)+2\varepsilon\alpha R(x,z,z,y)-\alpha^2K(y,z)\le (1-\alpha^2)c(p)
$$
for $|\alpha|<1$ and
$$
	K(x,z)+2\varepsilon\alpha R(x,z,z,y)-\alpha^2K(y,z)\ge (1-\alpha^2)c(p)
$$
for $|\alpha|>1$, where $\varepsilon = g(z,z)$. By continuity 
$$
	K(x,z)+2\varepsilon\alpha R(x,z,z,y)-\alpha^2K(y,z)=0
$$
holds good for \ $\alpha =\pm 1 $ \ and hence \ $R(x,z,z,y)=0$. So the Theorem follows
from Remark 4.

\vspace{0.1in}
R\,e\,m\,a\,r\,k 6. In Theorem 2 the requirement $K(\sigma)\le c(p)$ can be replaced by $K(\sigma)\ge c(p)$.

\vspace{0.1in}
R\,e\,m\,a\,r\,k 7. The conclusion of Theorem 2 rests true if  $|K(\sigma)|\le c(p)$ for any
antiholo\-morphic 2-plane $\sigma$ of signature $(+,-)$. Analogously for antiholomorphic
2-planes of signature $(+,+)$ or $(-,-)$ in the case of appropriate signature of $M$.

\vspace{0.7in}
\centerline{\large R E F E R E N C E S}

\vspace{0.1in}
\noindent
1. M. B\,a\,r\,r\,o\,s, A. R\,o\,m\,e\,r\,o. Indefinite Kaehlerian manifolds. {\it
Math. Ann}., {\bf 261}, 1982, 

55-62.

\noindent
2. M. D\,a\,j\,c\,z\,e\,r, K. N\,o\,m\,i\,z\,u. On sectional curvature of indefinite metrics.
{\it Math. Ann.},

 {\bf 247}, 1980, 279-282.

\noindent
3. M. D\,a\,j\,c\,z\,e\,r, K. N\,o\,m\,i\,z\,u. On the boundedness of Ricci curvature of an
indefinite 

metric.  {\it Bol. Soc. Brasil. Mat.}, {\bf 11}, 1980, 25-30.

\noindent
4. S. G. H\,a\,r\,r\,i\,s. A triangle comparison theorem for Lorentz manifolds. {\it Indiana 
Univ. 

Math. J.,} {\bf 31}, 1982, 289-308.
 
\noindent
5. R. S. K\,u\,l\,k\,a\,r\,n\,i. The values of sectional curvature of indefinite metrics.
{\it Comment. 

Math. Helvetici}, {\bf 54}, 1979, 173-176.

\vspace {0.3cm}
\noindent
{\it Center for mathematics and mechanics \ \ \ \ \ \ \ \ \ \ \ \ \ \ \ \ \ \ \ \ \ \ \ \ \ \ \ \ \ \ \ \ \ \
Received 01.03.1984

\noindent
1090 Sofia   \ \ \ \ \ \ \ \ \ \ \ \ \ \ \ \ \  P. O. Box 373}

\end{document}